\documentclass[psamsfonts]{amsart}

\usepackage{amscd,amsmath,amssymb,amsfonts,verbatim}
\usepackage{times}
\usepackage[cmtip, all]{xy}
\usepackage{graphicx}
\usepackage[active]{srcltx}
\usepackage{color}

\usepackage[color, notref, notcite]{showkeys}     % refs and labels
\definecolor{refkey}{gray}{.5}   % graylevel for refs
\definecolor{labelkey}{gray}{.5} % graylevel for labels
\definecolor{Red}{rgb}{1,0,0}

\newcommand{\pf}{{\bf Proof : }}

\newcommand{\qedwhite}{\hfill \ensuremath{\Box}}

\newtheorem{theo}{Theorem}[section]
\newtheorem{prop}[theo]{Proposition}
\newtheorem{lem}[theo]{Lemma}
\newtheorem{defi}[theo]{Definition}
\newtheorem{cor}[theo]{Corollary}

\newtheorem{notn}[theo]{Notation}
\newtheorem{rem}[theo]{Remark}
\newtheorem{con}[theo]{Convention}
\newtheorem{obs}[theo]{Observation}

\title{Homotopy and Commutativity Principle}
\author{Ravi A. Rao and Sampat Sharma}
\newcommand{\Addresses}{{% additional braces for segregating \footnotesize
  \bigskip
  \footnotesize
  
  \textsc{Ravi A. Rao, School of Mathematics, Tata Institute of Fundamental 
Research,\\  
           1, Dr. Homi Bhabha Road, Mumbai 400005, INDIA}\par\nopagebreak
  \textit{E-mail:} Ravi A. Rao \texttt{<ravi@math.tifr.res.in>}

  \medskip
  
 \textsc{Sampat Sharma, School of Mathematics, Tata Institute of Fundamental Research,\\ \noindent
           1, Dr. Homi Bhabha Road, Mumbai 400005, INDIA
   } \par\nopagebreak
  \textit{E-mail}: Sampat ~Sharma \texttt{<sampat@math.tifr.res.in>}

  \medskip

  }}

\begin{document}

\maketitle

\subjclass 2010 Mathematics Subject Classification:{13C10, 13H99, 19B10, 19B14.}
           
~ {Homotopy, Classical ~groups, ~Transvection ~groups.}

\section{\bf{Introduction}}

$R$ will denote a commutative ring with $1 \neq 0$ in this article, unless stated otherwise. 

The subject of injective stability for the linear group (i.e. $K_1(R)$) began 
in the famous paper of 
Bass--Milnor--Serre (\cite{bms}) where it was shown, in essence, that large 
sized stably elementary matrices were actually elementary matrices. This was 
shown by showing that the sequence (of 
pointed sets) 
\begin{eqnarray*}
\cdots \longrightarrow \frac{SL_n(R)}{E_n(R)} \longrightarrow \frac{SL_{n+1}(R)}{E_{n+1}(R)} 
\longrightarrow \cdots
\end{eqnarray*}
stabilizes.  The estimate they got was $n = 3$, when $\dim(R) = 1$, and 
for $n \geq {\rm max}\{3, d + 3\}$ otherwise. They conjectured that the 
correct bound for the linear quotients should be  
$n \geq {\rm max}\{3, d + 2\}$; which was established by L.N. Vaserstein 
in \cite{V}. 

In \cite{4} A.A. Suslin established the normality of the elementary 
linear subgroup $E_n(R)$ in $GL_n(R)$, for $n \geq 3$. This was a major 
surprise at that time as it was known due to the work of P.M. Cohn in 
\cite{cohn} that in general $E_2(R)$ is not normal in $GL_2(R)$. This is 
the initial precursor to study the non-stable $K_1$ groups $\frac{SL_n(R)}{E_n(R)},~n\geq 3$.

This theorem can also be got as a consequence of the Local-Global Principle of 
D. Quillen (for projective modules) in \cite{qu}; and its analogue for the 
linear group of elementary matrices $E_n(R[X])$, when $n \geq 3$ due to 
A. Suslin in \cite{4}. In fact, in \cite{22} it is shown that, in some sense,
the normality property of the elementary group $E_n(R)$ in $SL_n(R)$ is 
equivalent to having a Local-Global Principle for $E_n(R[X])$.  

In {\cite{1}}, A. Bak proved the following beautiful result:

\begin{theo} $($A. Bak$)$ For an almost commutative ring $R$ with identity with 
centre $C(R)$. The group $\frac{SL_{n}(R)}{E_{n}(R)}$ 
 is nilpotent of class atmost $\delta(C(R)) + 3 - n$, where  
$\delta(C(R)) < \infty ~and~ n\geq 3$, where  $\delta(C(R))$ is 
the Bass--Serre- dimension of $C(R)$.
\end{theo}

This theorem, which is proved by a localisation and completion technique, 
which evolved from an adaptation of the proof of the 
Suslin's $K_1$-analogue of Quillen's  
Local-Global Principle, was the starting point of our investigation. In this 
paper, we show (see Corollary \ref{2.20})

\begin{theo} 
\label{ab}
Let $R$ be a local ring, and let $A = R[X]$. Then the 
group $\frac{SL_{n}(A)}{E_{n}(A)}$ is an abelian group for $n \geq 3$. 
\end{theo} 

This theorem is a simple consequence of the following principle (see 
Theorem \ref{2.19}):

\vskip0.15in

\begin{theo} $($Homotopy and commutativity principle$):$ 
Let $R$ be a commutative ring. Let $\alpha \in SL_n(R)$, $n \geq 
3$, be homotopic to the identity. Then, for any $\beta \in SL_n(R)$, 
$\alpha \beta = \beta \alpha \varepsilon$, for some $\varepsilon \in E_n(R)$. 
\end{theo}

\vskip0.15in

This principle is a consequence of the Quillen--Suslin Local-Global principle; 
and using a {\it non-symmetric} application of it as done by A. Bak in 
\cite{1}.  

\vskip0.15in

The existence of a Local-Global Principle enables us to prove similar results 
in various groups. 

We restrict ourselves to the classical symplectic, 
orthogonal groups (and their relative versions); and to the automorphism groups
of a projective module (with 
a unimodular element), a symplectic module (with a hyperbolic summand), and an 
orthogonal module (with a hyperbolic symmand). 

However, our results can be 
extended to other Chevalley groups, relative Chevalley groups, reductive 
groups, etc. where such Local-Global Principles exist due to results of E. Abe
in \cite{abe}, A. Stepanov in \cite{as}, \cite{st}, Asok--Hoyois--Wendt in 
\cite{haw}, A. Stavrova in \cite{stv}, respectively. 

We could show that the 
symplectic quotients were abelian, but we could only establish that the 
orthogonal quotients are solvable of length atmost two. We do believe that the 
orthogonal quotient groups are also abelian; and prove this when the base 
ring is a regular local ring containing a field.

In (\cite{vdk2}, Theorem 4.1), W. van der Kallen has described an abelian 
group structure on the orbit space of unimodular rows under elementary 
action $\frac{Um_n(R)}{E_n(R)}$, when 
$n \geq 3$ and $d \leq 2n -4$, where $d$ is the dimension of $R$. In 
the paper \cite{vdk1} he does it in the case when $n = d+1$, where $d$ is 
the dimension of $R$; thereby extending the seminal work of L.N. Vaserstein 
in (\cite{7}, Theorem 5.2), when $d = 2$. His estimates come from similar 
estimates being true in case when $R$ is the ring of continuous real valued 
functions on a compact space $X$.

Let $Comp_r(R)$ denote the subset of $Um_r(R)$ consisting of the (completable) 
unimodular rows which can be completed to a matrix of determinant one. 
One of the interesting application of Theorem \ref{ab} is that the 
orbit set of completable unimodular rows over $R[X]$, when $R$ is a local 
ring, modulo the elementary action has an abelian group structure under 
matrix multiplication. (See Theorem \ref{2.25}.)

In particular, if one believes that the Bass--Suslin conjecture that unimodular
rows over a polynomial extension of a local ring is true, then one would have 
an abelian group structure on the orbit space $\frac{Um_n(R[X])}{E_n(R[X])}$. The only restriction on size is $n\geq 3$. Since 
one does know the truth of the Bass--Suslin conjecture when dimension $R$ is 
$3$ and $2R=R$ (see \cite{bq3}, \cite{bs3}); one does get $\frac{Um_3(R[X])}{E_3(R[X])}$ has an 
abelian group structure, when $R$ is a local ring of dimension $3$ in which 
$2$ is invertible.

Is there a (perhaps $\mathbb{A}^1-$homotopy) interpretation of this result 
from a topological point of view? 

\section{\bf{Linear and Symplectic group}}

First we collect some definitions and some known results,
and set notations which will be used throughout the paper.

\begin{defi}{\bf{Special linear group}}~{$SL_{n}(R)$}: The subgroup of the General
linear group $GL_{n}(R)$, of $n\times n$ invertible matrices of determinant
$1$.

\end{defi}

\begin{defi}{\bf{Elementary group}}~{$E_{n}(R) :$} The subgroup of all matrices 
of $GL_{n}(R)$ generated by $\{e_{ij}(\lambda) : \lambda \in R$,~ for~ $i\neq j\}$,
 where $e_{ij}(\lambda) = I_{n} +  \lambda E_{ij}$ and $E_{ij}$ is the 
matrix with $1$ on the $ij^{th}$ place and $0$'s elsewhere.
 
\end{defi}

\begin{notn}

Let $\psi_{1} = \begin{bmatrix}
                 0 & 1\\
                 -1 & 0\\
                \end{bmatrix},~~ \psi_{n} = \psi_{n-1} \perp \psi_{1};$~~ and 
                $\phi_{1} = \begin{bmatrix}
                 0 & 1\\
                  1 & 0\\
                \end{bmatrix},~~ \phi_{n} = \phi_{n-1} \perp \phi_{1}$, for $n > 1$.
\end{notn}

\begin{notn}
 Let $\sigma$ be the permutation of the natural numbers given by $\sigma(2i) = 2i-1$ and $\sigma(2i-1) = 2i$.
\end{notn}

\begin{defi}{\bf{Symplectic group}}~{$Sp_{2m}(R) :$} the group of all $2m \times 2m $ matrices 
$\{\alpha \in GL_{2m}(R)~\mid \alpha^{t}\psi_{m}\alpha = \psi_{m}\}$.
 
\end{defi}

\begin{defi}{\bf{Elementary Symplectic group}}~{$ESp_{2m}(R)$}: We define for $ 1\leq i \neq j\leq 2m,~z\in R,$\\
$$
se_{ij}(z)=
\begin{cases}
I_{2m} + zE_{ij},~~~\textit{if}~ i = \sigma(j);\\
 I_{2m} + zE_{ij} - (-1)^{i+j}zE_{\sigma(j)\sigma(i)}, ~~~ \textit{if}~ i\neq \sigma(j).

\end{cases}
$$

It is easy to verify that all these matrices belong to $Sp_{2m}(R)$. We call them the elementary symplectic matrices
over $R$. The subgroup generated by them is called the elementary symplectic group and is denoted by $ESp_{2m}(R)$.
 
\end{defi}

\begin{defi}{\bf{Orthogonal group}}~{$O_{2m}(R) :$} the group of all $2m \times 2m $ matrices 
$\{\alpha \in GL_{2m}(R)~\mid \alpha^{t}\phi_{m}\alpha = \phi_{m}\}$.
\end{defi}
\begin{defi}{\bf{Elementary Orthogonal group}}~{$EO_{2m}(R) :$} We define for $ 1\leq i \neq j\leq 2m,~z\in R,$\\
$$o_{ij}(z) = I_{2m} + zE_{ij} - zE_{\sigma(j)\sigma(i)}, ~~~\textit{if} ~ i\neq \sigma(j).$$\\
It is easy to verify that all these matrices belong to $O_{2m}(R)$. We call them the elementary orthogonal matrices over
$R$. The subgroup generated by them is called the elementary orthogonal group and is denoted by $EO_{2m}(R)$.
\end{defi}
\begin{notn}
 Let $R$ be a commutative ring with identity. In this paper $M(n,R)$ will denote the set of all $n\times n$ matrices over $R$,
 $G(n,R)$ will denote either 
 the linear group $GL_{n}(R)$ or the symplectic group $Sp_{2m}(R)$ , where $2m = n$. $E(n,R)$ will denote either elementary subgroups 
$E_{n}(R)$ or elementary symplectic subgroup $ESp_{2m}(R)$. And, $S(n,R)$ will denote either the special linear group 
 $SL_{n}(R)$ or the symplectic group $Sp_{2m}(R)$.
\end{notn}

\begin{con}
 Throughout this paper, we will assume size of the matrix is $n\geq 3$ in the linear case,  $n\geq 4$ in symplectic
 and $n\geq 6$ in orthogonal case, unless stated otherwise. 
\end{con}

\begin{lem}
\label{2.11}
 $($L.N. Vaserstein$)$ $(${\cite[Lemma 5.5]{7}}$)$ For an associative ring $R$ with identity, and for any natural number $m$
 $$E_{2m}(R)e_{1} = (Sp_{2m}(R) \cap E_{2m}(R))e_{1}.$$
\end{lem}

\begin{rem}
\label{2.12}
 It was observed in $(${\cite[Lemma 2.13]{3''}}$)$  that Vaserstein's proof actually shows that  
 $E_{2m}(R)e_{1} = ESp_{2m}(R)e_{1}.$
\end{rem}$~~~~~~~$ In view of above remark, or otherwise, one has:

\begin{lem}
\label{2.13}
$(${\cite[Chapter 1, Proposition 5.4]{2''}}$)$
Let $ c = (c_{1},\ldots, c_{n})$ be a unimodular row over a 
semilocal ring $R$. Then $(c_{1},\ldots, c_{n})\in e_{1} E(n,R); ~for~n\geq 2,$~ i.e.
$$(c_{1},\ldots, c_{n}) \overset{E(n, R)} {\sim}  (1,0,\ldots, 0)
 ~for ~n \geq 2.$$
\end{lem}

$~~~~~~~$ The next Lemma is well-known. We include it with a proof, for completeness.
\begin{lem}
\label{2.14}
  $($Only for the linear and the symplectic group$)$ Let $R$
  be a local ring. For $n\geq 2, S(n,R) = E(n,R)$, where $n = 2m,~ m$ is any
  natural number.
 \end{lem}
${\pf}$ For the linear case we prove the result by induction on $n$. When $n = 2$, it is obvious as 
$SL_{2}(R) = E_{2}(R)$. For $n>2$, 
let $\alpha \in SL_{n}(R)$,
By Lemma \ref{2.13},
$$ \alpha ~~ {\overset{E_{n}(R)}\sim} ~~\begin{bmatrix}
                 1 & 0\\
                 0 & \alpha^{'}\\
                \end{bmatrix}, ~~~~~~~~~\alpha^{'}\in E_{n}(R).$$ By induction hypothesis we have $\alpha^{'}\in E_{n-1}(R)$, 
                thus $\alpha \in E_{n}(R)$.
     \par
  ~~~~~~In the symplectic case let $\tau \in Sp_{2m}(R)$. We use the induction on $m$, for 
  $m=1$, $SL_{2}(R) = E_{2}(R) = Sp_{2}(R) = ESp_{2}(R)$.
  Since $\tau e_{1} \in E_{2m}(R)e_{1}$, by 
 Lemma \ref{2.11} (and remark following it), $\tau e_{1} \in ESp_{2m}(R)e_{1}$. Let $\tau e_{1} = \varepsilon_{1}e_{1}$ 
 for some 
   $\varepsilon_{1} \in ESp_{2m}(R)$.
    Hence $\varepsilon_{1}^{-1}\tau e_{1} = e_{1}.$ Therefore, we can find $\varepsilon_{2} \in ESp_{2m}(R)$ such that 
    $\varepsilon_{2}^{-1}\varepsilon_{1}^{-1}\tau = I_{2} \perp \tau^{\ast}$ for some $\tau^{\ast} \in Sp_{2m-2}(R)$. By induction $\tau^{\ast} \in ESp_{2m-2}(R)$. Repeating this process we can reduce $\tau$ to a 
     $2\times 2$ symplectic matrix.\\ 
    $~~~~~~~~~~~~~~~~~~~~~~~~~~~~~~~~~~~~~~~~~~~~~~~~~~~~~~~~~~~~~~~~~~~~~~~~~~~~~~~~~~~~~~~~~~~~~~~~~~~~~~~~~~~~~~\qedwhite$
         
 $~~~~~~~$ We begin with some initial observations:

\begin{lem}
\label{2.15}
 Let $R$ be a local ring and $\alpha(X),~\beta(X) \in S(n,R[X])$. Then the commutator, 
 $$[\alpha(X),~\beta(X)] \in [\alpha(X)\alpha(0)^{-1},~~\beta(X)\beta(0)^{-1}]E(n, ~R[X])$$
\end{lem}
${\pf}$ Since $R$ is a local ring, $S(n,R) = E(n,R)$ for all $n\geq 2$, by Lemma \ref{2.14}. Thus $\alpha(0), \beta(0) 
\in E(n,R)$. 
Let 
$ s = \alpha(X)\alpha(0)^{-1},~ t = \beta(X)\beta(0)^{-1}.$ Then, 
\begin{align*}
[\alpha(X),~\beta(X)] & =  [\alpha(X)\alpha(0)^{-1}\alpha(0),~~\beta(X)\beta(0)^{-1}\beta(0)]\\
 &= s\alpha(0)t\beta(0)(s\alpha(0))^{-1}(t\beta(0))^{-1}\\
& = 
sts^{-1}t^{-1}(tst^{-1}\alpha(0)ts^{-1}t^{-1})(ts\beta(0)\alpha(0)^{-1}s^{-1}t^{-1})
(t\beta(0)^{-1}t^{-1}).
\end{align*}

Since $E(n, R[X])$ is a normal subgroup of $S(n, R[X])$, 
hence $(tst^{-1}\alpha(0)ts^{-1}t^{-1})$,\\ 
$(ts\beta(0)\alpha(0)^{-1}s^{-1}t^{-1})$, $(t\beta(0)^{-1}t^{-1}) \in E(n, R[X]).$\\
$~~~~~~~~~~~~~~~~~~~~~~~~~~~~~~~~~~~~~~~~~~~~~~~~~~~~~~~~~~~~~~~~~~~~~~~~~~~~~~~~~~~~~~~~~~~~~~~~~~~~~~~~~~~~~~\qedwhite$

\begin{theo}
\label{2.16}
 {\bf(Local-Global Principle for the Linear groups)} $(${\cite[Theorem 3.1] {4}}$)$ Let $R$ be a commutative ring, 
 $n\geq 3$ and $\alpha \in GL_{n}(R[X])$ such that 
 $\alpha(0) = \textit{Id}$. Then $\alpha$ lies in $E_{n}(R[X])$ if and only if for every maximal ideal $\mathfrak{m}$ of $R$, the
  canonical image of $\alpha$ in $GL_{n}(R_{\mathfrak{m}}[X])$ lies in $E_{n}(R_{\mathfrak{m}}[X])$.
\end{theo}

\begin{theo}
\label{2.17}
{\bf (Local-Global Principle for the Symplectic groups)} 
$(${\cite[Theorem 3.6] {kop}}$)$
 Let $m\geq 2$ and $\alpha(X) \in Sp_{2m}(R[X])$, with $\alpha(0) = \textit{Id}$. Then $\alpha(X) \in ESp_{2m}(R[X])$ 
 if and only if for any maximal ideal $\mathfrak{m} \subset R$, the canonical image of $\alpha(X) \in 
 Sp_{2m}(R_{\mathfrak{m}}[X])$ lies in $ESp_{2m}(R_{\mathfrak{m}}[X])$.
\end{theo}

For a uniform proof of above Theorems see $(${\cite {22}}$)$
\begin{defi}
 Let R be a ring. A matrix $\alpha \in S(n,R)$ is said to be 
 homotopic to identity if there exists a matrix $\gamma(X) \in 
 S(n, R[X])$ such that $\gamma(0) = \textit{Id}~and~ \gamma(1) = \alpha$.
\end{defi}

\begin{theo}
\label{2.19}
 Let $\alpha \in S(n,R)$ be homotopic to identity. Then $[\alpha, \beta] \in E(n,R)$, $\forall ~\beta \in S(n,R)$. 
\end{theo}

${\pf}$  Since $\alpha$ is homotopic to identity, there exists $\gamma \in S(n, R[X])$ such that $\gamma(0) = \textit{Id}, ~ 
\gamma(1) = \alpha$.
 Define,  $$\delta(X) = [\gamma(X),~\beta].$$
 \par 
 Note that $\delta(0) = \textit{Id}$, and for every maximal ideal $\mathfrak{m}$ of $R$,
$$\delta(X)_{\mathfrak{m}} = [\gamma(X)_{\mathfrak{m}},~\beta_{\mathfrak{m}}].$$
\par
By  Lemma \ref{2.14}, $\beta_{\mathfrak{m}} \in E(n, R_{\mathfrak{m}})$ and since $E(n,R)$ is normal in $S(n,R)$, 
we have $\delta(X)_{\mathfrak{m}} \in E(n,R_{\mathfrak{m}}[X])$. Thus by Theorem \ref{2.16} (respectively Theorem \ref{2.17}), 
$\delta(X) = [\gamma(X), \beta]
\in E(n, R[X]),$ which implies 
$$\delta (1) = [\gamma(1), \beta] = [\alpha, \beta] \in E(n, R).$$
$~~~~~~~~~~~~~~~~~~~~~~~~~~~~~~~~~~~~~~~~~~~~~~~~~~~~~~~~~~~~~~~~~~~~~~~~~~~~~~~~~~~~~~~~~~~~~~~~~~~~~~~~~~~~~~\qedwhite$
\begin{cor}
\label{2.20}
 Let $R$ be a local ring. Then the group $\frac{S(n, R[X])}{E(n, R[X])}$ is an abelian group.
\end{cor}
${\pf}$ Let $\alpha(X),\beta(X) \in S(n, R[X])$, we need to prove that $[\alpha(X), \beta(X)] \in E(n, R[X]).$ 
In view of Lemma 
\ref{2.15}, we may assume that $\alpha(0) = \beta(0) = \textit{Id}$. 
\par 
Define,  $\gamma(X,T) = \alpha(XT).$ Clearly $\gamma(X, 0) = \textit{Id} ~ and~ \gamma(X, 1) = \alpha(X)$; thus $\alpha(X)$ is 
homotopic to identity.
 Thus, one gets the desired result by Theorem \ref{2.19}.\\
$~~~~~~~~~~~~~~~~~~~~~~~~~~~~~~~~~~~~~~~~~~~~~~~~~~~~~~~~~~~~~~~~~~~~~~~~~~~~~~~~~~~~~~~~~~~~~~~~~~~~~~~~~~~~~~\qedwhite$

\subsection*{The Relative case}
$~~~~~~~$ Let $I$ be an ideal of a ring $R$, we shall denote by $GL_{n}(R,I)$ the kernel of the canonical mapping 
$GL_{n}(R)\longrightarrow GL_{n}\left(\frac{R}{I}\right).$ Let $SL_{n}(R,I)$ denotes the subgroup of $GL_{n}(R,I)$ of elements of determinant $1$.

\begin{defi}${\bf{The~ Relative~ Groups~ E_{n}(I),~ E_{n}(R,I):}}$
 Let $I$ be an ideal of $R$. The elementary group $E_{n}(I)$ is the subgroup of $E_{n}(R)$ generated as a group by the elements
 $e_{ij}(x),~x\in I,~1\leq i\neq j\leq n.$\\
  The relative elementary group $E_{n}(R,I)$ is the normal closure of $E_{n}(I)$ in $E_{n}(R)$.
  
\end{defi}

\begin{defi}${\bf{The~ Relative~ Groups~ ESp_{2m}(I),~ ESp_{2m}(R,I):}}$
 Let $I$ be an ideal of $R$. The elementary symplectic group $ESp_{2m}(I)$ is the subgroup of $ESp_{2m}(R)$ 
 generated as a group by the 
 elements $se_{ij}(x),~x\in I,~1\leq i\neq j\leq 2m.$
 \par
  The relative elementary symplectic group $ESp_{2m}(R,I)$ is the normal closure of $ESp_{2m}(I)$ in $ESp_{2m}(R)$.
  \end{defi}
  
  \begin{defi}${\bf{The~ Relative~ Groups~ EO_{2m}(I),~ EO_{2m}(R,I):}}$
 Let $I$ be an ideal of $R$. The elementary orthogonal group $EO_{2m}(I)$ is the subgroup of $EO_{2m}(R)$ 
 generated as a group by the 
 elements $oe_{ij}(x),~x\in I,~1\leq i\neq j\leq 2m.$
 \par
  The relative elementary orthogonal group $EO_{2m}(R,I)$ is the normal closure of $EO_{2m}(I)$ in $EO_{2m}(R)$.
  \end{defi}
  
\begin{notn}
  Let $R$ be a commutative ring with identity and $I$ be an ideal of $R$. In this paper $E(n,R,I)$ will denote either 
  relative elementary group $E_{n}(R,I)$
   or relative elementary symplectic group $ESp_{2m}(R,I),$ and $S(n, R, I)$ will denote either $SL_{n}(R,I)$ or the 
   relative symplectic 
   group $Sp_{2m}(R,I)$ where $2m = n$.

\end{notn}

\begin{defi}${\bf {Excision ~Ring:}}$
 Let $R$ be a ring and $I$ be an ideal of $R$. The excision ring $R\oplus I$, has coordinate wise addition and multiplication
 is given as follows:
 $$(r,i).(s, j) = (rs, rj+si+ij),~where~r,s \in R ~and~i,j\in I.$$
 \par
 The multiplicative identity of this group is $(1,0)$ and the additive identity is $(0,0)$.
\end{defi}

\begin{lem}
\label{2.29}
$($Anjan Gupta$)$ $($see {\cite[Lemma 4.3]{ggr}}$)$
 Let $(R, \mathfrak{m})$ be a local ring. Then the excision ring $R\oplus I$ with respect to a proper ideal $I\subsetneq R$ is 
 also a local ring with maximal ideal $\mathfrak{m}\oplus I$.
\end{lem}

\begin{lem}
\label{2.30}
 Let $R$ be a local ring and $I$ be a proper ideal of $R$ (i.e. $I\neq R$). Then $S(n,R,I) = E(n,R,I)~for~all ~n\geq 1.$
\end{lem}
${\pf}$ Let $\sigma \in S(n,R,I)$, we can write $\sigma = \textit{Id} + \sigma^{'}$, for some $\sigma^{'} \in M_{n}(I)$.
Let $\overset{\sim}{\sigma} = (\textit{Id}, \sigma^{'}) \in S(n, R\oplus I, 0\oplus I)$. By Anjan's Lemma, $R\oplus I$ is a 
local ring, thus by Lemma 
\ref{2.14}, 
\begin{eqnarray*}
\overset{\sim}{\sigma} &\in& E(n, R\oplus I) \cap S(n, R\oplus I, 0\oplus I) = E(n, R\oplus I, 0\oplus I)
\end{eqnarray*}
as $\frac{R\oplus I}{0\oplus I} 
\simeq R$ is a retract of $R\oplus I$. Thus,\\
$$\overset{\sim}{\sigma} = \prod_{k=1}^{m}\beta_{k}ge_{i_{k}j_{k}}(0,a_{k})\beta_{k}^{-1},~~~~\beta_{k}\in 
E(n, R\oplus I),~ a_{k}\in I.$$
\par
Now, consider the homomorphism
\begin{eqnarray*}
f:R\oplus I\longrightarrow R\\
(r,i)\longmapsto r+i.
\end{eqnarray*}
\par
This $f$ induces a map\\
$$\overset{\sim}{f}:E(n, R\oplus I, 0\oplus I)\longrightarrow E(n,R)$$
\par
Clearly,\begin{align*} \sigma &= \overset{\sim}{f}(\overset{\sim}{\sigma})\\ 
~~~~~~~~~~~~~~~~~~~~~~~~~~~~~~~~~~~~~~~~~~~ &= \prod_{k=1}^{m}\gamma_{k}ge_{i_{k}j_{k}}(0 + a_{k})\gamma_{k}^{-1}\\
&= \prod_{k=1}^{m}\gamma_{k}ge_{i_{k}j_{k}}(a_{k})\gamma_{k}^{-1} \in E(n, R, I);~~\mbox{since} ~ a_{k}\in I,
\end{align*}

 $ \noindent \mbox{where},~\gamma_{k} = \overset{\sim}{f}(\beta_{k})$\\
 $~~~~~~~~~~~~~~~~~~~~~~~~~~~~~~~~~~~~~~~~~~~~~~~~~~~~~~~~~~~~~~~~~~~~~~~~~~~~~~~~~~~~~~~~~~~~~~~~~~~~~~~~~~~~~~\qedwhite$

 \begin{theo}
 \label{rel}
 
Let $R$ be a ring and I be a proper ideal (i.e. $I\neq R$) of $R$. Let $\alpha \in S(n, R, I)$
which is homotopic to identity relative to an extended ideal. 
Then $[\alpha, \beta] \in E(n, R, I), \forall \beta \in S(n, R, I)$.
\end{theo}

${\pf}$ Let $\alpha, \beta \in S(n, R, I)$ be such that $\alpha$ is homotopic to identity relative to an extended ideal. 
We can write $\alpha =
\textit{Id} + \alpha^{'}, \beta = \textit{Id} + \beta^{'}$ for some $\alpha^{'}, \beta^{'} \in  M_{n}(I).$ 
Let $\sigma = [\alpha, \beta] = \textit{Id} + \sigma^{'}$ for some $\sigma^{'} \in M_{n}(I)$.  
Let $\overset{\sim}{\sigma} = (\textit{Id}, \sigma^{'}) \in S(n, R\oplus I, 0\oplus I)$. In view of Theorem 2.19,
 \begin{eqnarray*}
\overset{\sim}{\sigma} &\in& E(n, R\oplus I) \cap S(n, R\oplus I, 0\oplus I) = E(n, R\oplus I, 0\oplus I)
\end{eqnarray*}
as $\frac{R\oplus I}{0\oplus I} 
\simeq R$ is a retract of $R\oplus I$. Thus,
$$\overset{\sim}{\sigma} = \prod_{k=1}^{m}\varepsilon_{k}ge_{i_{k}j_{k}}(0,a_{k})\varepsilon_{k}^{-1},~~~~\varepsilon_{k}\in 
E(n, R\oplus I),~ a_{k}\in I.$$
\par
Now, consider the homomorphism
\begin{eqnarray*}
f:R\oplus I\longrightarrow R\\
(r,i)\longmapsto r+i.
\end{eqnarray*}
\par
This $f$ induces a map\\
$$\overset{\sim}{f}:E(n, R\oplus I, 0\oplus I)\longrightarrow E(n,R)$$
\par
Clearly,\begin{align*} \sigma &= \overset{\sim}{f}(\overset{\sim}{\sigma})\\ 
~~~~~~~~~~~~~~~~~~~~~~~~~~~~~~~~~~~~~~~~~~~ &= \prod_{k=1}^{m}\gamma_{k}ge_{i_{k}j_{k}}(0 + a_{k})\gamma_{k}^{-1}\\
&= \prod_{k=1}^{m}\gamma_{k}ge_{i_{k}j_{k}}(a_{k})\gamma_{k}^{-1} \in E(n, R, I);~~\mbox{since} ~ a_{k}\in I,
\end{align*}

 $ \noindent \mbox{where},~\gamma_{k} = \overset{\sim}{f}(\varepsilon_{k})$\\
 $~~~~~~~~~~~~~~~~~~~~~~~~~~~~~~~~~~~~~~~~~~~~~~~~~~~~~~~~~~~~~~~~~~~~~~~~~~~~~~~~~~~~~~~~~~~~~~~~~~~~~~~~~~~~~~\qedwhite$ \\
$~~~~~~~~$ In view of the well-known Swan--Weibel homotopy trick $(${\cite [Appendix 3] {2''}}$)$, one has:
\begin{cor}
\label{grad}
  Let $A$ = $\bigoplus_{d\geq 0}A_{d}$ be a graded ring with augmentation ideal $A_{+} = \bigoplus_{d\geq 1}A_{d}$. 
Then $\frac{S(n, A,A_{+})}{E(n,A,A_{+})}$ is an abelian 
 group.
\end{cor}

${\pf}$ Consider the ring homomorphism 
$$\varphi : A\longrightarrow A[T]$$ 
$$~~~~~a_{0} + a_{1}+\cdots \longmapsto a_{0}+ a_{1}T+\cdots$$
\par
Note that $\varphi$ is an injective ring homomorphism. For any element $\alpha = (\alpha_{ij}) \in S(n,A, A_{+})$, 
define $\alpha(T) = (\varphi(\alpha_{ij}))$. Now, note that $\alpha(0) = \textit{Id}$ and $\alpha(1) = \alpha$. Thus by 
Theorem \ref{rel} $\frac{S(n, A ,A_{+})}{E(n, A,A_{+})}$ is an abelian group.\\
$~~~~~~~~~~~~~~~~~~~~~~~~~~~~~~~~~~~~~~~~~~~~~~~~~~~~~~~~~~~~~~~~~~~~~~~~~~~~~~~~~~~~~~~~~~~~~~~~~~~~~~~~~~~~~~\qedwhite$
\begin{cor}
 Let $A$ be an affine algebra of dimension $d\geq 2$ over a perfect $C_{1}$ field $k$ and $(d+1)!$ is a unit in $k$
 . Let $\sigma \in SL_{d+1}(A)$ be a stably elementary 
 matrix. Then $[\sigma, \tau] \in E_{d+1}(A)$, for all $\tau \in SL_{d+1}(A).$
\end{cor}

${\pf}$ In view of $(${\cite [Theorem 3.4] {3'}}$)$, there exists a matrix $\sigma(X) \in SL_{d+1}(A[X])$ with $\sigma(0) = 
\textit{Id}$ and
 $\sigma(1) = \sigma.$ Now, we are through by Theorem \ref{2.19}.\\
 $~~~~~~~~~~~~~~~~~~~~~~~~~~~~~~~~~~~~~~~~~~~~~~~~~~~~~~~~~~~~~~~~~~~~~~~~~~~~~~~~~~~~~~~~~~~~~~~~~~~~~~~~~~~~~~\qedwhite$
\begin{cor}
  Let $A$ be an affine algebra of dimension $d\geq 3$ over an algebrically closed field $k$ and ${d!}$ is a unit in $k$.
 Then the group $\frac{SL_{d}(A)\cap E_{d+1}(A)}{E_{d}(A)}$ is an abelian group.
\end{cor}
${\pf}$ Let $\sigma \in {SL_{d}(A)\cap E_{d+1}(A)}$. In view of $(${\cite [Corollary 7.7] {frs}}$)$, $\sigma$ is 
homotopic to identity. Thus we are through by Theorem \ref{2.19}.\\
$~~~~~~~~~~~~~~~~~~~~~~~~~~~~~~~~~~~~~~~~~~~~~~~~~~~~~~~~~~~~~~~~~~~~~~~~~~~~~~~~~~~~~~~~~~~~~~~~~~~~~~~~~~~~~~\qedwhite$

\begin{cor}
 Let $A$ be an affine algebra of even dimension $d$ over a field $k$ of cohomological dimension $\leq 1$. 
 If $(d+1)!A = A$ and $4|d$, then $\frac{Sp_{d}(A)\cap ESp_{d+2}(A)}{ESp_{d}(A)}$ is an 
 abelian group.
\end{cor}
${\pf}$ Let $\sigma \in {Sp_{d}(A)\cap ESp_{d+2}(A)}$. In view of  $(${\cite [Theorem 1]{brp}}$)$, 
$\sigma$ is symplectic homotopic to $\textit{Id}$. Thus we are through by Theorem \ref{2.19}.\\
$~~~~~~~~~~~~~~~~~~~~~~~~~~~~~~~~~~~~~~~~~~~~~~~~~~~~~~~~~~~~~~~~~~~~~~~~~~~~~~~~~~~~~~~~~~~~~~~~~~~~~~~~~~~~~~\qedwhite$

The reader should contrast the next result with the results $(${\cite [Theorem 5.1]{3'}}$)$ and $(${\cite [Proposition 7.10]{vdk2}}$)$ of W. van der Kallen.

\begin{theo}\label{2.25}
  Let $A$ = $\bigoplus_{d\geq 0}A_{d}$ be a graded ring with augmentation ideal $A_{+} = \bigoplus_{d\geq 1}A_{d}$. 
Then for $n\geq 3$, $\frac{Comp_{n}(A,A_{+})}{E_{n}(A,A_{+})}$ has an abelian 
 group structure under matrix multiplication. In particular, for $n\geq 3$, the first row map
 $${SL_{n}(A, A_{+})}\longrightarrow \frac{Comp_{n}(A, A_{+})}{E_{n}(A, A_{+})}$$
 $$\sigma\longmapsto[e_{1}\sigma]$$
 is a group homomorphism.
\end{theo}

${\pf}$ Since $\frac{SL_{n}(A, A_{+})}{E_{n}(A, A_{+})}$ is an abelian group, it is enough to prove that matrix
 multiplication gives a well defined (abelian) operation on $\frac{Comp_{n}(A, A_{+})}{E_{n}(A, A_{+})}$.
  Let $v, u \in \frac{Comp_{n}(A, A_{+})}{E_{n}(A, A_{+})}$ such that
$$v = e_{1}{\alpha} = e_{1}{\alpha^{'}};~~~~~\alpha, \alpha^{'}\in SL_{n}(A, A_{+})$$
$$u = e_{1}{\beta} = e_{1}{\beta^{'}};~~~~~\beta, \beta^{'}\in SL_{n}(A, A_{+})$$
\par
To get a well-defined multiplication on $\frac{Comp_{n}(A, A_{+})}{E_{n}(A, A_{+})}$, we need 
to prove that $[e_{1}\alpha \beta] = [e_{1}\alpha^{'}\beta^{'}]$. By Corollary \ref{grad}, 
$$[e_{1}\alpha \beta] = [e_{1}\alpha^{'}\beta] = [e_{1}\beta \alpha^{'}] = [e_{1}\beta{'}\alpha^{'}] = 
[e_{1}\alpha^{'}\beta^{'}].$$
$~~~~~~~~~~~~~~~~~~~~~~~~~~~~~~~~~~~~~~~~~~~~~~~~~~~~~~~~~~~~~~~~~~~~~~~~~~~~~~~~~~~~~~~~~~~~~~~~~~~~~~~~~~~~~~\qedwhite$

\begin{cor}
 Let $A$ be a commutative ring and $Comp_{n}(A)$ denote the subset of $Um_{n}(A)$ consisting of those unimodular rows which 
 can be completed 
 to an invertible matrix of determinant $1$.  If $R$ is a local ring, then $\frac{Comp_{n}(R[X])}{E_{n}(R[X])}$ has an abelian group structure ,
 under matrix multiplication for $n\geq 3.$
\end{cor}

\begin{rem}
 There exist examples which show that $\frac{Comp_n(A, A_+)}{E_n(A, A_+)}$ is non-trivial, for 
some graded ring $A$:
\par
Let $R = k[X, Y, Z]/(Z^7 - X^2 - Y^3)$, where $k$ is $\mathbb{C}$ or any sufficiently 
large field of characteristic $\neq 2$.
It is shown in $($\cite{fsr}, page 4$)$ that 
if $B = R[T, T^{-1}]$, then there is a maximal ideal $\mathfrak{m}$ 
for which $NW_E(B_\mathfrak{m}) \neq 0$. By $($\cite{SV}, Theorem 5.2 (b)$)$ the 
Vaserstein symbol $\frac{Um_3(B_\mathfrak{m}[W])}{E_3(B_\mathfrak{m}[W])} 
\longrightarrow W_E(B_\mathfrak{m}[W])$ is onto. Hence, there exists a unimodular row $v(W) \in Um_3(B_\mathfrak{m}[W])$ 
which is not elementarily completable. However, by \cite{bs3}, $v(W)$ is 
completable. 
\end{rem}

\begin{theo} \label{2.31}
$(${\bf Local~Global~ Principle~ for ~Extended ~Ideals}$)$ $(${\cite [Theorem 1.3]
{2'''}}$)$
 Let $\alpha(X) \in G(n, R[x],I[X])$ be such that $\alpha(0) = \textit{Id}$. If 
 $\alpha_{\mathfrak{m}}(X) \in E(n, R_{\mathfrak{m}}[X], I_{\mathfrak{m}}[X])$ for 
 every maximal $\mathfrak{m}$ of $R$, then $\alpha(X) \in E(n, R[X], I[X]).$
\end{theo}
\begin{lem}
\label{2.32}
 Let $R$ be a local ring and $I$ be a proper ideal of $R$ (i.e. $I\neq R$). Then the group $\frac{S(n, R[X], I[X])}{E(n, R[X], I[X])}$ is an abelian group.
\end{lem}
${\pf}$ Let $\alpha(X),\beta(X) \in S(n,R[X], I[X])$, we need to prove that $[\alpha(X), \beta(X)]
\in E(n,R[X],I[X]).$ In view of Lemma 
\ref{2.15}, we may assume that $\alpha(0) = \beta(0) = \textit{Id}$. Define, 
 $$\gamma(X,T) = [\alpha(XT), \beta(X)]$$
 \par
 Note that, $\gamma(X,0) = \textit{Id}$ and for every maximal ideal $\mathfrak{m}$ of $R[X]$, 
$\gamma(X,T)_{\mathfrak{m}} = [\alpha(XT)_{\mathfrak{m}}, \beta(X)_{\mathfrak{m}}]$, by Lemma \ref{2.30}, 
$\beta(X)_{\mathfrak{m}} \in 
E(n, R[X]_{\mathfrak{m}}, I[X]_{\mathfrak{m}})
\subseteq E(n, R[X]_{\mathfrak{m}}, I[X]_{\mathfrak{m}}[T]);$ and since $E_{n}$ is normal in $SL_{n}$, we have,
$$\gamma(X,T)_{\mathfrak{m}} \in E(n, R[X]_{\mathfrak{m}}, I[X]_{\mathfrak{m}}[T]).$$
\par
Thus by  Theorem \ref{2.31}, $\gamma(X,T) \in E((n,R[X], I[X])[T])$ which implies that
$\gamma(X,1) = [\alpha(X), \beta(X)] \in E(n, R[X], I[X]).$\\
$~~~~~~~~~~~~~~~~~~~~~~~~~~~~~~~~~~~~~~~~~~~~~~~~~~~~~~~~~~~~~~~~~~~~~~~~~~~~~~~~~~~~~~~~~~~~~~~~~~~~~~~~~~~~~~\qedwhite$

\begin{cor}
 Let $A$ be an affine algebra of dimension $d$ over an algebrically closed field $k$ and $I = (a)$ be a principal ideal. 
 Assume $(d+1) ! \in k^{\ast}$, $d \equiv 1~(\textit{mod}~ 4)$. Then $\frac{Sp_{d-1}(A,I)\cap ESp_{d+1}(A,I)}{ESp_{d-1}(A,I)}$ is an 
 abelian group.
\end{cor}
${\pf}$ Let $\sigma \in {Sp_{d-1}(A,I)\cap ESp_{d+1}(A,I)}$. In view of $(${\cite [Theorem 5.4]{anjan}}$)$, 
$\sigma$ is symplectic homotopic to $\textit{Id}$. Thus we are through by Theorem \ref{2.19}.\\
$~~~~~~~~~~~~~~~~~~~~~~~~~~~~~~~~~~~~~~~~~~~~~~~~~~~~~~~~~~~~~~~~~~~~~~~~~~~~~~~~~~~~~~~~~~~~~~~~~~~~~~~~~~~~~~\qedwhite$

\section{\bf{Transvection groups}}
\par 
First, we collect some definitions and some known results and set notations which will be used in this paper.

\begin{defi}
 Let $M$ be a finitely generated module over a ring $R$. An element $m$ of $M$ is said to be unimodular in $M$ if
 $Rm \cong R $ and $ M \cong Rm \oplus M^{'}$,
 for some $R$-submodule $M^{'}$ of $M$.
\end{defi}

\begin{defi}
 For an element $m \in M$, one can attach an ideal, called the order ideal of $m$ in $M$, viz. 
 $$O_{M}(m) = \{f(m)\mid f\in M^{\ast} =Hom(M,R)\}$$
\end{defi}

\begin{defi}
We define a transvection of a finitely generated $R$-module as follows: Let $M$ be a finitely generated R-module. 
Let $q\in M $ and 
$f\in M^{\ast}$ with $f(q) = 0$. An automorphism of $M$ of the form $1 + f_{q}$ (defined by $f_{q}(p) = f(p)q,$ for $p \in M$), 
will be
 called a {\bf{transvection}} of $M$ if either $q\in Um(M)$ or $f\in Um(M^{\ast}).$ We denote by 
 $Trans(M)$ the subgroup of $Aut(M)$
 generated by 
 transvections of $M$.
\end{defi}

\begin{defi}
  Let $M$ be a finitely generated $R$-module . The automorphisms of the form $(p,a) \mapsto (p + ax, a)$ and $(p,a) \mapsto 
  (p, a + f(p)),$ where 
 $x\in M $ and $f\in M^{\ast}$, are called {\bf elementary transvections} of $M\oplus R$. $($Note that we can regard $f$ as 
 an element of $(M\oplus R)^{\ast}$ by defining $f(0,1) = 0$.$)$ By taking $q = (x,0)$ and $f\in (M\oplus R)^{\ast}$ such that 
 $f: (y,t)  \mapsto t$ for
  $(y,t) \in (M\oplus R),$ one can verify that the automorphism $(p,a)\mapsto (p+ax, a)$ is in $Trans(M\oplus R)$. 
  Similarly, by taking $
  q = (0,1)~and ~ f\in (M\oplus R)^{\ast}$ such that $f : (0,1)\mapsto 0$ one can verify that the automorphism 
  $(p,a)\mapsto (p, a+f(p))$ 
  is in $Trans(M\oplus R).$
  \par 
  The subgroup of $Trans(M\oplus R)$ generated by the elementary transvections is denoted by $ETrans(M\oplus R)$.

\end{defi}

\begin{defi}
 A {\bf symplectic (respectively orthogonal)} $R$-module is a pair $(P, \langle,\rangle)$, where $P$ is a projective $R$-module of even rank
 and $\langle,\rangle :
 P \times P\longrightarrow R$ is a non-degenerate alternating (respectively symmetric) bilinear form.
\end{defi}
\begin{defi}
  Let $(P_{1}, \langle,\rangle_{1})$ and $(P_{2}, \langle,\rangle_{2})$ be two symplectic (respectively orthogonal) $R$-modules. 
  Their {\bf orthogonal sum} is a pair
  $(P, \langle,\rangle)$, where $ P = P_{1} \oplus P_{2}$ and the inner product is defined by 
  $\langle (p_{1}, p_{2}), (q_{1}, q_{2})\rangle = 
  \langle p_{1}, q_{1}\rangle_{1} + \langle p_{2}, q_{2}\rangle_{2}.$ Since this form is also non-singular
  we shall henceforth denote 
  $(P, \langle,\rangle)$ by $P_{1} \perp P_{2}$ called the orthogonal sum 
  of $(P_{1}, \langle,\rangle_{1})$ and $(P_{2}, \langle,\rangle_{2})$ .
  
\end{defi}
\begin{defi}
 For a projective $R$-module $P$ of rank $n$, we define $\mathbb{H}(P)$ of rank $2n$ supported by $P \oplus P^{\ast}$, 
 with form 
 $\langle (p, f), (p^{'}, f^{'})\rangle  = f(p^{'}) - f^{'}(p)$  for the symplectic modules and $f(p^{'}) + f^{'}(p)$ for the
 orthogonal modules. 
\end{defi}
\begin{defi}
 An {\bf isometry} of a symplectic (respectively orthogonal) module $(P, \langle , \rangle )$ is an automorphism of $P$ 
 which fixes the bilinear form. The
 group of isometries of $(P, \langle , \rangle )$ is denoted by $Sp(P)$ for the symplectic modules and 
 $O(P)$ for the orthogonal modules.
\end{defi}

\begin{defi}
We define a symplectic transvection as follows: Let $\Psi : P\longrightarrow P^{\ast}$ be an induced isomorphism. 
Let $\alpha : R\longrightarrow P$
be a $R$-linear map defined by $\alpha(1) = u.$ Then $\alpha^{\ast}\Psi$ defined by $\alpha^{*}\Psi(p) = \langle u,p\rangle$ 
is in $P^{\ast}$.
Let $v\in P$ be such that $\alpha^{*}\Psi(v) = \langle u,v\rangle = 0$.
An automorphism $\sigma_{(u,v)}$ of $(P,\langle , \rangle)$ of the form
$$\sigma_{(u,v)}(p) = p + \langle u, p\rangle v +  \langle v, p\rangle u + \langle u, p\rangle u$$
for $u,v\in P$ with $\langle u, v\rangle = 0$ will be called a {\bf symplectic transvection} of $(P,\langle , \rangle)$ if 
either $v\in Um(P)$
or $\alpha^{\ast}\Psi \in Um(P^{\ast}).$\
   Since $\langle \sigma_{(u,v)}(p(1)), \sigma_{(u,v)}(p(2))\rangle = \langle p_{1}, p_{2}\rangle, \sigma_{(u,v)}\in 
   Sp(P,\langle, \rangle.$
 Note that $\sigma^{-1}_{(u,v)}(p) = p - \langle u, p\rangle v -  \langle v, p\rangle u - \langle u, p\rangle u$.\
 The subgroup of $Sp(P, \langle, \rangle)$ generated by symplectic transvections is denoted by $Trans_{Sp}(P)$.
\end{defi}

\begin{defi}
 The symplectic transvections of $P\perp R^{2}$ of the form
 $$(p,b,a)\longmapsto (p+aq, b-\langle p, q\rangle + a, a),$$
 $$(p,b,a)\longmapsto (p+bq, b, a - b +\langle p, q\rangle),$$
 where $a,b \in R$ and $p,q \in P$, are called {\bf elementary symplectic transvections}.
 \par
 One can verify that above two maps belong to $Trans_{Sp}(P\perp R^{2}).$ The subgroup of $Trans_{Sp}(P\perp R^{2})$
 generated by elementary symplectic transvections is denoted by $ETrans_{Sp}(P\perp R^{2}).$
 \par
  In a similar manner we can find a transvection $\tau_{(u,v)}$ for an orthogonal module $(P, \langle , \rangle )$. 
  For this we 
  need to assume that $u,v\in P$ are {\bf isotropic}, i.e. $\langle u,u\rangle = \langle v,v\rangle = 0.$
  \end{defi}
  
\begin{defi}
 An automorphism $\tau_{(u,v)}$ of $(P, \langle , \rangle )$ of the form\
 $$ \tau_{(u,v)}(p) = p - \langle u,p\rangle v + \langle v,p\rangle u$$ for $u, v \in P~with~ 
 \langle u,v\rangle = \langle u,u\rangle =
 \langle v,v\rangle = 0$ will be called an {\bf isotropic orthogonal transvection} of $(P, \langle , \rangle )$ if either 
 $v\in Um(P)$ or $\alpha^{\ast}\Psi \in Um(P^{\ast})$.
 \par 
 One can verify that $\tau_{(u,v)} \in O(P, \langle , \rangle )$ 
  and $ \tau^{-1}_{(u,v)}(p) = p + \langle u,p\rangle v - \langle v,p\rangle u$.
  The subgroup of $O(P, \langle , \rangle )$ generated by  isotropic orthogonal transvections is denoted by $Trans_{O}(P)$.
 \end{defi}

\begin{defi}
 The isotropic orthogonal transvections of $(P \perp R^{2})$ of the form \
 $$(p,b,a)\longmapsto (p-aq, b+\langle p, q\rangle , a),$$
 $$(p,b,a)\longmapsto (p-bq, b, a - \langle p, q\rangle),$$
 where $a, b \in R$ and $p,q\in P$, are called {\bf elementary orthogonal transvections}.
 \par 
 The subgroup of $Trans_{O}(P\perp R^{2})$ generated by the elementary orthogonal 
  transvections is denoted by $ETrans_{O}(P\perp R^{2})$.

\end{defi}
\begin{notn}
\label{3.14}
 In this paper $P$ will denote either a finitely generated projective module of rank $n$, a symplectic module or
  an orthogonal module of even rank $n = 2m$
 with a fixed form $\langle, \rangle$. And Q will denote $P\oplus R$ in the 
 linear case and $P\perp R^{2}$ otherwise. We assume that $n \geq 2$, when dealing
 with linear case and symplectic case and $n\geq 4$ otherwise. We use notation $G(Q)$ to
  denote $Aut(Q)$, $Sp(Q, \langle, \rangle)$  respectively; $S(Q)$ will denote $SL(Q) = 
  \{\sigma \in Aut(Q): ~~\wedge^{n}\sigma = 1\}$, 
  $Sp(Q, \langle, \rangle)$ respectively; $T(Q)$ to denote $Trans(Q)$, $Trans_{Sp}(Q)$ 
   respectively; and $ET(Q)$ to denote $ETrans(Q)$, $ETrans_{Sp}(Q)$ respectively.
\end{notn}

\begin{theo}
\label{3.15}
 {\bf(Local-Global Principle for Transvection Groups)} $(${\cite [Theorem 3.6] {2}}$)$ Let $R$ be a commutative ring with 
 identity and $Q$ be as in Notation \ref{3.14}.
 Suppose $\sigma(X) \in G(Q[X])$ with $\sigma(0) = \textit{Id}.$ If for every maximal ideal 
 $\mathfrak{m}$ of $R$, 
 $$
\sigma_{\mathfrak{m}}(X) \in
\begin{cases}
E(n+1,~~R_{\mathfrak{m}}[X]) ~for~linear~case,~\\
 
 E(n+2,~~R_{\mathfrak{m}}[X]) ~ otherwise.
\end{cases}
$$
Then $\sigma(X) \in ET(Q[X])$.
\end{theo}

\begin{theo}
$(${\cite [Theorem 2] {2}}$)$
 $T(Q) = ET(Q)$. Hence $ET(Q)$ is normal subgroup of $G(Q)$.
\end{theo}

\begin{theo}
\label{3.16}
 Let $\sigma \in S(Q)$ such that $\sigma$ is homotopic to identity. Then $[\sigma, \tau] \in ET(Q)$ for all $\tau \in S(Q)$.
\end{theo}
${\pf}$ Since $\sigma$ is homotopic to identity there exists $\varphi(X) \in S(Q[X])$ such that 
$\varphi(0) = \textit{Id}~and ~\varphi(1) = \sigma$. Define
$$\varPsi(X) = [\varphi(X), \tau]$$
\par 
Note that $\varPsi(0) = \textit{Id} $ and for every maximal ideal $\mathfrak{m}$ of $R$, $\varPsi(X)_{\mathfrak{m}} = 
[\varphi(X)_{\mathfrak{m}}, \tau_{\mathfrak{m}}]$. By Lemma \ref{2.14}, $\sigma_{\mathfrak{m}} \in E(n+1, R_{\mathfrak{m}})$ in 
linear case and 
$\sigma_{\mathfrak{m}} \in E(n+2, R_{\mathfrak{m}})$ in symplectic case. Since $E(n,R)$ is normal in $S(n,R)$
, we have $\varPsi(X)_{\mathfrak{m}} \in
 E(n+1, R_{\mathfrak{m}}[X])$ in linear case and $\varPsi(X)_{\mathfrak{m}} \in
 E(n+2, R_{\mathfrak{m}}[X])$ in symplectic case. Thus by Theorem \ref{3.15}, $\varPsi(X) \in ET(Q[X])$ which implies
 $$\varPsi(1) = [\varphi(1), \tau] = [\sigma, \tau] \in ET(Q).$$
 $~~~~~~~~~~~~~~~~~~~~~~~~~~~~~~~~~~~~~~~~~~~~~~~~~~~~~~~~~~~~~~~~~~~~~~~~~~~~~~~~~~~~~~~~~~~~~~~~~~~~~~~~~~~~~~\qedwhite$

\begin{cor}
 Let $R$ be a commutative ring and $P$ be a finitely generated projective $R$-module of rank $n = 2m$. 
 Then the group $\frac{S(Q[X], (X))}{ET(Q[X], (X))}$ is an abelian group.
\end{cor}
${\pf}$  Let $\sigma(X),~ \tau(X) \in S(Q[X])$, we need to prove that $[\sigma(X),~ \tau(X)] \in ET(Q[X])$.
Define, $$\gamma(X) = [\sigma(X), \tau(X)]$$
 \par
 Then $\gamma(0) = \textit{Id}$. For every maximal ideal $\mathfrak{m}$ of $R$, by Corollary \ref{2.20},
$[\sigma(X)_{\mathfrak{m}}, \tau(X)_{\mathfrak{m}}] \in  E(n+1, R_{\mathfrak{m}}[X])$ in linear case 
and $[\sigma(X)_{\mathfrak{m}}, \tau(X)_{\mathfrak{m}}] \in E(n+2, R_{\mathfrak{m}}[X])$ in symplectic case. 
Hence by  Theorem \ref{3.15}, $\gamma(X) \in ET(Q[X], (X))$.\\
$~~~~~~~~~~~~~~~~~~~~~~~~~~~~~~~~~~~~~~~~~~~~~~~~~~~~~~~~~~~~~~~~~~~~~~~~~~~~~~~~~~~~~~~~~~~~~~~~~~~~~~~~~~~~~~\qedwhite$

\begin{theo}
\label{3.18}
 {\bf(Local-Global Principle for Transvection Groups in Relative case)} $(${\cite [Theorem 1.3] {2'''}}$)$ Let $R$ be a 
 commutative ring with identity and $Q$ be as in Notation \ref{3.14}.
 Suppose $\sigma(X) \in G(Q[X], I[X])$ with $\sigma(0) = \textit{Id}.$ If for every maximal ideal 
 $\mathfrak{m}$ of $R$, 
 $$
\sigma_{\mathfrak{m}}(X) \in
\begin{cases}
E(n+1,~~R_{\mathfrak{m}}[X],I_{\mathfrak{m}}[X]) ~for~linear~case,~\\
 
 E(n+2,~~R_{\mathfrak{m}}[X],I_{\mathfrak{m}}[X]) ~ otherwise.
\end{cases}
$$
 Then $\sigma(X) \in ETrans(Q[X], I[X])$.
 \end{theo}
  
\begin{lem}
  \label{3.19}
 Let $R$ be a commutative ring and $I$ be a proper ideal of $R$ (i.e. $I\neq R$) and $P$ be a finitely generated projective 
 $R$-module of rank $n = 2m$. 
 Then the group $\frac{S(Q[X], XI[X])}{ET(Q[X], XI[X])}$ is an abelian group.
\end{lem}
 ${\pf}$  Let $\sigma(X),\tau(X) \in S(Q[X], XI[X])$, we need to prove that $[\sigma(X), \tau(X)] \in ET(Q[X],XI[X]).$ 
 Define, $$\gamma(X) = [\sigma(X), \tau(X)]$$
 \par
 Note that, $\gamma(0) = \textit{Id}$. For every maximal ideal $\mathfrak{m}$ of $R$,   by Lemma \ref{2.32},
$[\sigma(X)_{\mathfrak{m}}, \tau(X)_{\mathfrak{m}}] \in \\ E(n+1, R_{\mathfrak{m}}[X], XI_{\mathfrak{m}}[X])$ in linear case 
and $[\sigma(X)_{\mathfrak{m}}, \tau(X)_{\mathfrak{m}}] \in E(n+2, R_{\mathfrak{m}}[X], XI_{\mathfrak{m}}[X])$ in symplectic case. 
Thus $\gamma(X) \in ET(Q[X], XI[X])$ by Theorem \ref{3.18}.\\
$~~~~~~~~~~~~~~~~~~~~~~~~~~~~~~~~~~~~~~~~~~~~~~~~~~~~~~~~~~~~~~~~~~~~~~~~~~~~~~~~~~~~~~~~~~~~~~~~~~~~~~~~~~~~~~\qedwhite$

\section{\bf{Orthogonal groups}}
\par
  Throughout this section we will assume that $1/2 \in R$, where $R$ is a commutative ring with identity.

\begin{defi}{\bf{Lower Central Series}}
 Let $G$ be a group, and define $G_{0} = G$, $G_{n} = [G_{n-1}, G]~for~ n\geq 1.$ With these notations, we have
 $$G = G_{0}\supseteq G_{1} \supseteq \cdots G_{n} \supseteq \cdots.$$
 The above series of subgroups of $G$ is called the lower central series of group $G$.
\end{defi}
\par 
We say a group $G$ is nilpotent if lower central series terminates after finitely many terms and if $G_{n}$ 
is the first subgroup which is trivial in the series then $G$ is said to be nilpotent of nilpotency class $n$. 

\begin{defi}{\bf{Derived Series}}
 Let $G$ be a group, and define $G^{0} = G$, $G^{n} = [G^{n-1}, G^{n-1}]~for~ n\geq 1.$ With these notations, we have
 $$G = G^{0}\supseteq G^{1} \supseteq \cdots G^{n} \supseteq \cdots.$$
 The above series of subgroups of $G$ is called the derived series of group $G$.
\end{defi}
\par 
We say a group $G$ is solvable if derived series terminates after finitely many terms and if $G^{n}$ 
is the first subgroup which is trivial in the series then $G$ is said to be a solvable group of length $n$. 

\begin{defi}
 By $EO_{R}(Q \perp \mathbb{H}(P))~\cdotp O_{R}(\mathbb{H}(P))$ we shall mean the subset 
 $\{ \sigma_{1}\sigma_{2}\lvert \sigma_{1}\in EO_{R}(Q \perp \mathbb{H}(P))
 , \sigma_{2} \in O_{R}(\mathbb{H}(P))\}$ of $O_{R}(Q \perp \mathbb{H}(P))$.
\end{defi}
\begin{theo}
\label{4.3}
 {\bf(Local-Global Principle for the Orthogonal groups)} $(${\cite[Theorem 4.2] {sus}}$)$
 Let $m\geq 3$ and $\alpha(X) \in SO_{2m}(R[X])$, with $\alpha(0) = \textit{Id}$. Then $\alpha(X) \in EO_{2m}(R[X])$ 
 if and only if for any 
 maximal ideal $\mathfrak{m} \subset R$, the canonical image of $\alpha(X)$ in $SO_{2m}(R_{\mathfrak{m}}[X])$ 
 lies in $EO_{2m}(R_{\mathfrak{m}}[X])$.
\end{theo}

\begin{lem}
\label{4.4}
{$($R.A.Rao$)$}$ (${\cite [Lemma 2.2] {3}}$)$
 Let $R$ be a ring with Jacobson dimension $\leq d$. Let $(Q,q)$ be a diagonalisable quadratic $R$-space. Consider the
 quadratic $R$-space $Q \perp \mathbb{H}(P)$, where rank $P > d$. Then,
 $$O_{R}(Q \perp \mathbb{H}(P)) = 
 EO_{R}(Q \perp \mathbb{H}(P))~\cdotp O_{R}(\mathbb{H}(P)) = O_{R}(\mathbb{H}(P))~\cdotp EO_{R}(Q \perp \mathbb{H}(P))$$
\end{lem}
\begin{defi}{\bf{Spinor Norm}}
  Suppose $R$ be a local ring and $M$ be an $R$-module and $B$ be a non-degenerate symmetric bilinear form. Let $G$
  be the orthogonal group 
 corresponding to $B$. The spinor norm is a group homomorphism
 $$SN : G\longrightarrow \frac{(R^{\ast})}{(R^{\ast})^{2}}.$$
 \par
 The homomorphism is defined as follows: any element of $G$ arising as reflection orthogonal to vector $v$ is sent to the value 
 $B(v,v)$ modulo $(R^{\ast})^{2}$. This extends to a well-define and unique homomorphism on all of $G$. The reflection 
 orthogonal to 
 vector $v$ is defined as 
 \begin{align*}& \tau_{v}:  M\longrightarrow M\\
 & x\mapsto x-2v\frac{B(v,x)}{B(v,v)}, ~for~all ~v\in M.
 \end{align*}
\end{defi}

\begin{obs}
\label{4.6}
 One can write the matrix $\begin{bmatrix}
                         u & 0\\
                         0 & u^{-1}\\
                        \end{bmatrix}$ as a product $\tau_{e_{1}-e_{2}} \tau_{e_{1}-ue_{2}}$. Hence its spinor norm
                        is $4u$.
\end{obs}
\par 
In view of $(${\cite[Theorem 4] {kli}}$)$, $EO_{2m}(R)$ is a normal subgroup of $SO_{2m}(R)$ when
         $R$ is a local ring.
         
\begin{lem}
\label{4.7}
Let $R$ be a local ring and $I$ be a proper ideal of $R$ (i.e. $I\neq R$). 
If $\alpha \in SO_{2m}(R,I)$, then $\alpha^{2} \in EO_{2m}(R,I)$.
\end{lem}
${\pf}$ Let $\alpha \in SO_{2m}(R,I)$. Since $\alpha \in SO_{2m}(R,I)$,  
we can write $\alpha = \textit{Id} + \alpha'$, for some $\alpha' \in M_{2m}(I)$.
Let $\overset{\sim}{\alpha} = (\textit{Id}, \alpha') \in SO_{2m}(R\oplus I, 0\oplus I)$. 
By Lemma \ref{2.29}, $R\oplus I$ is a local ring. In view of Lemma \ref{4.4} $SO_{2m}(R\oplus I) = 
SO_{2}(R\oplus I) \cdotp EO_{2m}(R\oplus I)$. Since $R$ is a commutative ring and  by $(${\cite [Lemma 3.2] {dser}}$)$
every element of
$SO_{2}(R\oplus I)$ looks like 
$$\begin{bmatrix}
                         u & 0\\
                         0 & u^{-1}\\
                        \end{bmatrix}, ~~~~u\in (R\oplus I)^{\ast}.$$ 
\par Thus $\alpha^{2}$ looks like 
$$\begin{bmatrix}
                         u^{2} & 0\\
                         0 & u^{-2}\\
                        \end{bmatrix}, ~~~~u\in (R\oplus I)^{\ast}.$$ 
 \par In view of observation \ref{4.6}, or otherwise, spinor norm of $\alpha^{2}$ is $4u^{2}$, 
 a square in $(R\oplus I)^{\ast}$. 
    Thus by $(${\cite [Theorem 6] {klie}}$)$
     we have, $\alpha^{2} \in EO_{2}(R\oplus I)$. (The details of the proof can 
be found in \cite{kli}.)
\begin{eqnarray*}
\overset{\sim}{\alpha}^{2} &\in& EO_{2m}(R\oplus I) \cap SO_{2m}(R\oplus I, 0\oplus I) = EO_{2m}(R\oplus I, 0\oplus I)
\end{eqnarray*}
as $\frac{R\oplus I}{0\oplus I} 
\simeq R$ is a retract of $R\oplus I$. Thus,
$$\overset{\sim}{\alpha}^{2} = \prod_{k=1}^{t}b_{k}ge_{i_{k}j_{k}}(0,a_{k})b_{k}^{-1},~~~~b_{k}\in EO_{2m}(
R\oplus I),~ a_{k}\in I.$$
 Now, consider the homomorphism
$$f:R\oplus I\longrightarrow R$$
$$~~~~~~~~~~~~~(r,i)\longmapsto r+i.$$
This $f$ induces a map
$$\overset{\sim}{f}:EO_{2m}(R\oplus I, 0\oplus I)\longrightarrow EO_{2m}(R, I).$$ Thus $ \alpha^{2} \in EO_{2m}(R, I).$\\
$~~~~~~~~~~~~~~~~~~~~~~~~~~~~~~~~~~~~~~~~~~~~~~~~~~~~~~~~~~~~~~~~~~~~~~~~~~~~~~~~~~~~~~~~~~~~~~~~~~~~~~~~~~~~~~\qedwhite$

\begin{cor}
 Let $R$ be a local ring and $I$ be a proper ideal of $R$ (i.e. $I\neq R$). Then the group $\frac{SO_{2m}(R, I)}
 {EO_{2m}(R, I)}$ is an abelian
 group for all $m\geq 1$. In fact, every element of this group is of order $2.$
\end{cor}

\begin{lem}
\label{4.17}
 Let $R$ be a local ring then the group $\frac{SO_{2m}(R)}{EO_{2m}(R)}$ is an abelian group for all $m\geq 1$. In fact, 
 every element of this group 
 is of order $2$. In particular, $[SO_{2m}(R), SO_{2m}(R)] = EO_{2m}(R)$.
\end{lem}
${\pf}$ In view of R.A. Rao's Lemma $SO_{2m}(R) = SO_{2}(R) \cdotp EO_{2m}(R)$. Every element of
$SO_{2}(R)$ looks like 
$$\begin{bmatrix}
                         u & 0\\
                         0 & u^{-1}\\
                        \end{bmatrix}, ~~~~u\in R^{\ast}.$$ 
\par Thus $SO_{2}(R)$ is an abelian group which implies that $\frac{SO_{2m}(R)}{EO_{2m}(R)}$ is also an abelian group. 
Since, for every element 
$\alpha \in SO_{2}(R)$, $\alpha^{2}$ looks like 
$$\begin{bmatrix}
                         u^{2} & 0\\
                         0 & u^{-2}\\
                        \end{bmatrix}, ~~~~u\in R^{\ast}.$$ 
    \par Again spinor norm of $\alpha^{2}$ is $4u^{2}$, a square in $R^{\ast}$. 
    Thus by $(${\cite [Theorem 6] {klie}}$)$
     we have, $\alpha^{2} \in EO_{2}(R)$. (The details of the proof can 
be found in \cite{kli}.)\\
$~~~~~~~~~~~~~~~~~~~~~~~~~~~~~~~~~~~~~~~~~~~~~~~~~~~~~~~~~~~~~~~~~~~~~~~~~~~~~~~~~~~~~~~~~~~~~~~~~~~~~~~~~~~~~~\qedwhite$

\begin{theo}
\label{4.8}
 Let $R$ be a local ring. Then the group $\frac{SO_{2m}(R[X])}{EO_{2m}(R[X])}$ is
 a solvable group of length at most $2$.
\end{theo}
${\pf}$ Let $\alpha(X), \beta(X) \in [SO_{2m}(R[X]), SO_{2m}(R[X])]$, 
we need to prove $[\alpha(X), \beta(X)] \in EO_{2m}(R[X]).$ In view of 
Lemma \ref{4.17}, we may assume that $\alpha(0) = \beta(0) = \textit{Id}$.
Define, 
$$\gamma(X,T) = [\alpha(XT), \beta(X)].$$
\par For every maximal ideal $\mathfrak{m}$ of $R[X]$,
$$\gamma(X,T)_{\mathfrak{m}} = [\alpha(XT)_{\mathfrak{m}}, \beta(X)_{\mathfrak{m}}].$$
\par Since $\beta(X)_{\mathfrak{m}}\in [SO_{2m}(R[X]_{\mathfrak{m}}), SO_{2m}(R[X]_{\mathfrak{m}})] = EO_{2m}(R[X]_{\mathfrak{m}})$
 and  $EO_{2m}(R[X]_{\mathfrak{m}}) \trianglelefteq  SO_{2m}(R[X]_{\mathfrak{m}}),$ 
 thus
 $\gamma(X,T)_{\mathfrak{m}} \in EO_{2m}(R[X]_{\mathfrak{m}}[T])$ and $\gamma(X,0) = \textit{Id}$. Thus by Theorem \ref{4.3},
 $\gamma(X, T) \in EO_{2m}(R[X,T]),$ by putting $T = 1$, one gets $\gamma(X,1)  = [\alpha(X), \beta(X)] \in EO_{2m}(R[X]).$\\
 $~~~~~~~~~~~~~~~~~~~~~~~~~~~~~~~~~~~~~~~~~~~~~~~~~~~~~~~~~~~~~~~~~~~~~~~~~~~~~~~~~~~~~~~~~~~~~~~~~~~~~~~~~~~~~~\qedwhite$

\begin{theo}
\label{4.15}
 Let $R$ be a local ring and $I$ be a proper ideal of $R$ (i.e. $I\neq R$). Then the group
 $\frac{SO_{2m}(R[X], I[X])}{EO_{2m}(R[X], I[X])}$ is
 a solvable group of length at most $2$.
\end{theo}
${\pf}$ Let $\alpha(X), \beta(X) \in [SO_{2m}(R[X], I[X]), SO_{2m}(R[X], I[X])] $, 
we need to prove that $[\alpha(X), \beta(X)] \in EO_{2m}(R[X], I[X]).$ In view of 
Lemma \ref{4.8}, we may assume that $\alpha(0) = \beta(0) = \textit{Id}.$
Define,
$$\gamma(X,T) = [\alpha(XT), \beta(X)].$$
\par For every maximal ideal $\mathfrak{m}$ of $R[X]$,
$$\gamma(X,T)_{\mathfrak{m}} = [\alpha(XT)_{\mathfrak{m}}, \beta(X)_{\mathfrak{m}}].$$
\par Since $\beta(X)_{\mathfrak{m}}\in [SO_{2m}(R[X]_{\mathfrak{m}},I[X]_{\mathfrak{m}}), 
SO_{2m}(R[X]_{\mathfrak{m}},I[X]_{\mathfrak{m}})] = EO_{2m}(R[X]_{\mathfrak{m}},I[X]_{\mathfrak{m}}).$
By $(${\cite [Corollary 2.13] {sus}}$)$ $EO_{2m}(R[X]_{\mathfrak{m}}, I[X]_{\mathfrak{m}})$ is normal in
$SO_{2m}(R[X]_{\mathfrak{m}},I[X]_{\mathfrak{m}}),$ 
 and so
 $\gamma(X,T)_{\mathfrak{m}} \in EO_{2m}(R[X]_{\mathfrak{m}}[T],\\I[X]_{\mathfrak{m}}[T])$. Since $\gamma(X,0) = \textit{Id}$, 
  by Theorem \ref{2.31},
 $\gamma(X, T) \in EO_{2m}(R[X,T], I[X,T]),$ by putting $T = 1$, one gets $\gamma(X,1)  = [\alpha(X), \beta(X)] \in 
 EO_{2m}(R[X], I[X]).$\\
 $~~~~~~~~~~~~~~~~~~~~~~~~~~~~~~~~~~~~~~~~~~~~~~~~~~~~~~~~~~~~~~~~~~~~~~~~~~~~~~~~~~~~~~~~~~~~~~~~~~~~~~~~~~~~~~\qedwhite$

\begin{theo}
\label{4.9}
 Let $R$ be a commutative ring and $P$ be a finitely generated projective $R$-module of rank $n = 2m$. 
 Then the group $\frac{SO (Q[X], (X))}{ETrans_{O}(Q[X], (X))}$ is a solvable group of 
length at most $2$, for $m\geq 2$, where
 $Q = P\perp R^{2}$.
\end{theo}
${\pf}$ Let $\sigma(X), \tau(X) \in [SO(Q[X], (X)),~ SO(Q[X], (X))]$,
we need to prove that 
$[\sigma(X), \tau(X)] \in ETrans_{O}(Q[X], (X)).$ 
Define,
$$\gamma(X) = [\sigma(X),~\tau(X)]$$
\par
Then $\gamma(0) = \textit{Id}$. For every maximal ideal $\mathfrak{m}$ of $R$, $\sigma(X)_{\mathfrak{m}} 
\in [SO_{2m+2}(R_{\mathfrak{m}}[X]), ~
SO_{2m+2}(R_{\mathfrak{m}}[X])],
~ \tau(X)_{\mathfrak{m}} \in [SO_{2m+2}(R_{\mathfrak{m}}[X]), ~SO_{2m+2}(R_{\mathfrak{m}}[X])]$. In view of 
Theorem \ref{4.15}, $\gamma(X)_{\mathfrak{m}} \in 
EO_{2m+2}(R_{\mathfrak{m}}[X], (X)_{\mathfrak{m}}).$ Now using
Theorem \ref{3.15}, we get $\gamma(X) \in ETrans_{O}(Q[X], (X)).$\\
$~~~~~~~~~~~~~~~~~~~~~~~~~~~~~~~~~~~~~~~~~~~~~~~~~~~~~~~~~~~~~~~~~~~~~~~~~~~~~~~~~~~~~~~~~~~~~~~~~~~~~~~~~~~~~~\qedwhite$

\begin{cor}
 Let $R$ be a commutative ring and $I$ be a proper ideal of $R$ (i.e. $I\neq R$) and $P$ be a finitely generated
 projective $R$
 -module of rank $n = 2m$. 
 Then the group $\frac{SO (Q[X], XI[X])}{ETrans_{O}(Q[X], XI[X])}$ is a solvable group of length at most $2$
 , for $m\geq 2$, where $Q = P\perp R^{2}.$
\end{cor}
${\pf}$ Let $\sigma(X), \tau(X) \in [SO(Q[X], XI[X]),~ SO(Q[X], XI[X])]$, 
we need to prove that $[\sigma(X), \tau(X)] \\ \in ETrans_{O}(Q[X], XI[X]).$ 
Define,
$$\gamma(X) =  [\sigma(X),~\tau(X)]$$
\par
Then $\gamma(0) = \textit{Id}$. For every maximal ideal $\mathfrak{m}$ of $R$, $\sigma(X)_{\mathfrak{m}} \in
[SO_{2m+2}(R_{\mathfrak{m}}[X], XI_{\mathfrak{m}}[X]), ~SO_{2m+2}(R_{\mathfrak{m}}[X],\\ XI_{\mathfrak{m}}[X])],
~ \tau(X)_{\mathfrak{m}} \in [SO_{2m+2}(R_{\mathfrak{m}}[X], XI_{\mathfrak{m}}[X]), ~SO_{2m+2}(R_{\mathfrak{m}}[X],
XI_{\mathfrak{m}}[X])]$. In view of Theorem \ref{4.15}, $\gamma(X)_{\mathfrak{m}} \in 
EO_{2m+2}(R_{\mathfrak{m}}[X], XI_{\mathfrak{m}}[X]).$  Now using
Theorem \ref{3.18} we get $\gamma(X) \in ETrans_{O}(Q[X], XI[X]).$\\
$~~~~~~~~~~~~~~~~~~~~~~~~~~~~~~~~~~~~~~~~~~~~~~~~~~~~~~~~~~~~~~~~~~~~~~~~~~~~~~~~~~~~~~~~~~~~~~~~~~~~~~~~~~~~~~\qedwhite$

\begin{lem}
\label{4.16}
Let $R$ be a local ring. If $\alpha(X), \beta(X) \in SO_{2m}(R[X])$ with $\alpha(0) = \textit{Id}$, then 
$[\alpha(X), \beta(X)^{2}] \in EO_{2m}(R[X])$. 
\end{lem}
${\pf}$ Define,
$$\gamma(X,T) = [\alpha(XT), \beta(X)^{2}].$$
\par For every maximal ideal $\mathfrak{m}$ of $R[X]$,
$$\gamma(X,T)_{\mathfrak{m}} = [\alpha(XT)_{\mathfrak{m}}, \beta(X)_{\mathfrak{m}}^{2}].$$
\par In view of $(${\cite [Theorem 6] {klie}}$)$ $\beta(X)_{\mathfrak{m}}^{2}\in  EO_{2m}(R[X]_{\mathfrak{m}})$
 and  $EO_{2m}(R[X]_{\mathfrak{m}})$ is normal in $SO_{2m}(R[X]_{\mathfrak{m}}),$ 
 thus
 $\gamma(X,T)_{\mathfrak{m}} \in EO_{2m}(R[X]_{\mathfrak{m}}[T])$ and $\gamma(X,0) = \textit{Id}$. Thus by Theorem \ref{4.3},
 $\gamma(X, T) \in EO_{2m}(R[X,T]),$ by putting $T = 1$, one gets $\gamma(X,1)  = [\alpha(X), \beta(X)^{2}] \in EO_{2m}(R[X]).$\\
 $~~~~~~~~~~~~~~~~~~~~~~~~~~~~~~~~~~~~~~~~~~~~~~~~~~~~~~~~~~~~~~~~~~~~~~~~~~~~~~~~~~~~~~~~~~~~~~~~~~~~~~~~~~~~~~\qedwhite$
 
\begin{cor}
\label{4.22}
 Let $R$ be a local ring. Then the group $\frac{SO_{2m}(R[X], (X))}{EO_{2m}(R[X], (X))}$ is
 a nilpotent group of class at most $2$.
\end{cor}
${\pf}$ Let $\alpha(X) \in SO_{2m}(R[X], (X)) ~and ~\beta(X) \in [SO_{2m}(R[X], (X)), SO_{2m}(R[X], (X))]$, 
we need to prove that $[\alpha(X), \beta(X)] \in EO_{2m}(R[X], (X)).$ Since in any group $G$, a commutator 
$[x,y] = xyx^{-1}y^{-1} = (xyx^{-1})^{2}x^{2}(x^{-1}y^{-1})^{2} $ is 
a product of squares. Therefore $\beta(X)$ can be 
written as a product of squares; and the result follows from Lemma \ref{4.16}. 
\begin{flushleft} We also give an alternative proof for this Corollary:
\end{flushleft}

Let $\alpha(X) \in SO_{2m}(R[X], (X)) ~and ~\beta(X) \in [SO_{2m}(R[X], (X)), SO_{2m}(R[X], (X))]$, 
we need to prove that $[\alpha(X), \beta(X)] \in EO_{2m}(R[X], (X)).$
Define,
$$\gamma(X,T) = [\alpha(XT), \beta(X)].$$
\par For every maximal ideal $\mathfrak{m}$ of $R[X]$,
$$\gamma(X,T)_{\mathfrak{m}} = [\alpha(XT)_{\mathfrak{m}}, \beta(X)_{\mathfrak{m}}].$$
\par Since $\beta(X)_{\mathfrak{m}}\in [SO_{2m}(R[X]_{\mathfrak{m}}), SO_{2m}(R[X]_{\mathfrak{m}})] = 
EO_{2m}(R[X]_{\mathfrak{m}})$
 and  $EO_{2m}(R[X]_{\mathfrak{m}})\trianglelefteq  SO_{2m}(R[X]_{\mathfrak{m}}),$ 
 thus
 $\gamma(X,T)_{\mathfrak{m}} \in EO_{2m}(R[X]_{\mathfrak{m}}[T])$ and $\gamma(X,0) = \textit{Id}$ since $\alpha(0) =
 \textit{Id}$. Thus by Theorem \ref{4.3},
 $\gamma(X, T) \in EO_{2m}(R[X,T]),$ by putting $T = 1$, one gets $\gamma(X,1)  = [\alpha(X), \beta(X)] \in EO_{2m}(R[X]).$ 
 Since $\gamma(0,1) = \textit{Id}$, thus 
$\gamma(X,1)  = [\alpha(X), \beta(X)] \in EO_{2m}(R[X], (X)).$\\
 $~~~~~~~~~~~~~~~~~~~~~~~~~~~~~~~~~~~~~~~~~~~~~~~~~~~~~~~~~~~~~~~~~~~~~~~~~~~~~~~~~~~~~~~~~~~~~~~~~~~~~~~~~~~~~~\qedwhite$

\par
 We believe that the orthogonal quotients (in Theorem \ref{4.8} and Theorem \ref{4.9}) are abelian groups, 
 we show this when the base ring 
 is a regular local ring containing a field.

\begin{defi}
 Let $k$ be a field. A ring $R$ is said to be essentially of finite type over $k$ if $R = S^{-1}C$, with $S$
 is a multiplicatively 
 closed subset of $C$, and $C = k[x_{1},\ldots ,x_{m}]/I$ is a quotient ring of a polynomial ring over $k$.
\end{defi}

\begin{prop}
\label{4.12}
 Let $R$ be a smooth affine algebra over a field $k$. If $\alpha(X) \in SO_{2m}(R[X])$ with $\alpha(0) = \textit{Id}$, then 
 $\alpha(X) \in EO_{2m}(R[X])$, for $m\geq 3$.
\end{prop}
${\pf}$ Let $\gamma(X,T) = \alpha(XT) \in SO_{2m}(R[X,T])$, then $\gamma(X,0) = \textit{Id}$. Thus by homotopy invariance 
$($See {\cite  {jhs}}, {\cite [Corollary 1.12] {jh}}, {\cite [Theorem 9.8] {msl}}.$)$ we have $\gamma(X,1) = \alpha(X) \in 
EO_{2m}(R[X])$. (The reader may also consult \cite{dr} for a version which is 
suitable for this application.) (There is a version for reductive groups 
in {\cite{stv}} which may be of independent interest.) \\
$~~~~~~~~~~~~~~~~~~~~~~~~~~~~~~~~~~~~~~~~~~~~~~~~~~~~~~~~~~~~~~~~~~~~~~~~~~~~~~~~~~~~~~~~~~~~~~~~~~~~~~~~~~~~~~\qedwhite$

The following Lemma is well known: 
\begin{lem}
\label{4.13}
 Let $R$ ba an affine algebra over a field $k$. Suppose $R_{\mathfrak{p}}$ is regular local ring for some $\mathfrak{p} 
 \in Spec(R)$. Then there 
 exists $s \notin \mathfrak{p}$ such that $R_{s}$ is a regular ring.
\end{lem}
${\pf}$ Let $J$ be the Jacobian ideal $R$. Then $V(J) = Sing(R)$. Since $R_{\mathfrak{p}}$ is regular local ring, 
$J\nsubseteq \mathfrak{p}$. Choose 
$s\in J\backslash \mathfrak{p}$. 
Now for every $\mathfrak{q} \in Spec(R)$ with $s\notin \mathfrak{q}$, we have $J\nsubseteq \mathfrak{q}$. Since every prime 
ideal of $R_{s}$ looks like $\mathfrak{q_{s}}$, for some $\mathfrak{q} \in Spec(R)$ with $s\notin \mathfrak{q}$, 
we get $(R_{s})_{\mathfrak{q_{s}}} = 
R_{\mathfrak{q}}$ is a regular local ring. Hence $R_{s}$ is a regular ring.\\
$~~~~~~~~~~~~~~~~~~~~~~~~~~~~~~~~~~~~~~~~~~~~~~~~~~~~~~~~~~~~~~~~~~~~~~~~~~~~~~~~~~~~~~~~~~~~~~~~~~~~~~~~~~~~~~\qedwhite$

\begin{theo}
\label{4.14}
 Let $R$ be a regular local ring essentially of finite type over a field $k$. 
 If $\sigma(X) \in SO_{2m}(R[X])$, with $\sigma(0) = \textit{Id}$, 
  then $\sigma(X) \in EO_{2m}(R[X])$, for $m\geq 3$.
 \end{theo}
${\pf}$ In view of Lemma \ref{4.13}, for every $\mathfrak{p} \in Spec(R)$ there exists $s \in R\backslash \mathfrak{p}$ 
such that $R_{s}$ is a smooth algebra.
 Therefore, by Proposition \ref{4.12}, $\sigma_{s}(X) \in EO_{2m}(R_{s}[X])$, which implies 
 $\sigma_{\mathfrak{p}}(X) \in EO_{2m}(R_{\mathfrak{p}}[X]).$ Now, by Theorem \ref{4.3} we have 
  $\sigma(X) \in EO_{2m}(R[X])$.\\
$~~~~~~~~~~~~~~~~~~~~~~~~~~~~~~~~~~~~~~~~~~~~~~~~~~~~~~~~~~~~~~~~~~~~~~~~~~~~~~~~~~~~~~~~~~~~~~~~~~~~~~~~~~~~~~\qedwhite$

\begin{cor}
 Let $R$ be a geometric regular local ring containing a field. Then the group $\frac{SO_{2m}(R[X])}{EO_{2m}(R[X])}$ 
 is an abelian group for 
 $m\geq 3$.
\end{cor}

${\pf}$  Let $\alpha(X), \beta(X) \in SO_{2m}(R[X])$, we need to prove that $[\alpha(X), \beta(X)] \in EO_{2m}(R[X]).$
Define,
$$\gamma(X) = [\alpha(0), \beta(0)]^{-1}[\alpha(X), \beta(X)].$$
 \par We will proceed by induction on dim $R$. If $\dim R = 0$, then $R$ is a field and the result follows
from Proposition \ref{4.12}. Therefore 
we assume $\dim R\geq 1$. In {\cite {pop}}, D. Popescu showed that if $R$ is a geometric
regular local ring then it is a filtered inductive limit 
of regular local rings essentially of finite type over a field.
\par
 Clearly, $\gamma(0) = \textit{Id}$. Hence by Theorem \ref{4.14},
 $\gamma(X) \in EO_{2m}(R[X])$. Now, using Lemma \ref{4.17}, we get 
 $[\alpha(X), \beta(X)] \in EO_{2m}(R[X])$.\\
$~~~~~~~~~~~~~~~~~~~~~~~~~~~~~~~~~~~~~~~~~~~~~~~~~~~~~~~~~~~~~~~~~~~~~~~~~~~~~~~~~~~~~~~~~~~~~~~~~~~~~~~~~~~~~~\qedwhite$

\medskip
\noindent
{\bf Acknowledgement:} We thank the referee for his suggestions, indicating 
a proof of Corollary \ref{grad}, and for a quick review.

\Addresses

 \end{document}